\documentclass[11pt, twoside]{article}
 \usepackage{amsfonts, amssymb, amsmath, amsthm, latexsym, amscd}
\newcommand{\braket}[2]{\langle #1,#2 \rangle}

\newcommand{\Ga}{\Gamma}
\newcommand{\ga}{\gamma}

\newcommand{\del}{\delta}
\newcommand{\Del}{\Delta}

\newcommand{\tht}{\theta}

\newcommand{\f}{\frac}

\newcommand{\lo}{\longrightarrow}

\newcommand{\BR}{\mathbb{R}}

\newcommand{\pa}{\partial}

\theoremstyle{plain}
\newtheorem{Th}{Theorem}[section]
\newtheorem{Lem}[Th]{Lemma}
\newtheorem{Prop}[Th]{Proposition}

\theoremstyle{definition}

\begin{document}
\addtolength{\textheight}{0 cm} \addtolength{\hoffset}{0 cm}
\addtolength{\textwidth}{0 cm} \addtolength{\voffset}{0 cm}

\title{On a class of periodic   quasilinear Schr\"{o}dinger
equations involving critical  growth   in ${\BR}^2$}

\author{Abbas Moameni  \textsc \bf \footnote{Research is supported by a
Postdoctoral Fellowship at the University of British Columbia.}}
\maketitle
\begin{center}
{\small Department of Mathematics \\
\small  University of British Columbia\\
\small Vancouver, B.C., Canada \\
\small {\tt moameni@math.ubc.ca }}
\end{center}

\begin{abstract}

We consider the equation $- \Delta u+V(x)u-
k(\Del(|u|^{2}))u=g(x,u),  u>0,  x \in
  {\BR}^2,$ where $V:{\BR}^2\rightarrow
{\BR}$  and $g:{\BR}^2 \times {\BR}\rightarrow {\BR}$ are two
continuous $1-$periodic functions. Also, we assume $g$ behaves like
$\exp (\beta |u|^4)$ as $|u|\rightarrow \infty.$ We prove the
existence of at least one weak solution $u \in H^1({\BR}^2)$ with
$u^2 \in H^1({\BR}^2).$
 Mountain pass in a suitable Orlicz space together with Moser-Trudinger
 are  employed to establish this result. Such
equations arise when one seeks for standing wave solutions for the
corresponding quasilinear Schr\"{o}dinger equations. Schr\"{o}dinger
equations of this type have been studied as models of several
physical phenomena. The nonlinearity here corresponds to the
superfluid film equation in plasma physics.
\end{abstract}
\section{Introduction}

 We are concerned with the existence of positive
solutions for quasilinear elliptic equations in the entire space,
\begin{eqnarray*}
- \Delta u+V(x)u- k(\Del(|u|^{2}))u=g(x,u), \quad u>0, x \in
  {\BR}^2,
\end{eqnarray*}
where  $V:{\BR}^2\rightarrow [0,\infty)$ and $g:{\BR}^2 \times
{\BR}\rightarrow [0,\infty)$ are nonnegative continuous functions.
Solutions of this equation are related to the existence of standing
wave solutions for quasilinear Schr\"{o}dinger equations of the form
\begin{equation}
i\pa_t z=- \Del z+W(x)z-l(x,|z|)z-k\Del h(|z|^2)h'(|z|^2)z, \quad
x\in {\BR}^N,  N \geq 2,
\end{equation}
where $W(x)$  is a given potential, $k$ is a real constant and $l$
and $h$ are real functions. Quasilinear equations of the form (1)
have been established  in several areas of physics corresponding to
various types of $h$.
 The  superfluid film equation in
plasma physics has this structure   for $h(s)=s$, ( Kurihura in
[8]). In the case $h(s)=(1+s)^{1/2}$, equation (1) models the
self-channeling of a high-power ultra short laser in matter, see
 [21]. Equation (1) also
appears in  fluid mechanics [8,9], in the theory of Heidelberg
ferromagnetism and magnus [10], in dissipative quantum mechanics and
in condensed matter theory [14]. We consider the case  $h(s)=s$ and
$k>0$. Setting $z(t,x)=\exp(-iFt)u(x)$ one obtains a corresponding
equation of elliptic type which has the formal variational
structure:
\begin{eqnarray}
- \Delta u+V(x)u-k(\Del(|u|^{2}))u=g(x,u), \quad u>0, x \in
  {\BR}^N,
\end{eqnarray}
where $V(x)=W(x)-F$ is the new potential function and $g$ is the new nonlinearity.\\

 Note that, for the case $g(u)=|u|^{p-1}u$ with $N\geq 3,$
$p+1=2(2^*)= \frac{4N}{N-2}$
 behaves like a critical exponent for the above equation [13, Remark
 3.13]. For the
 subcritical case $p+1 <22^*$ the existence of solutions for problem (2) was studied in [10, 11, 12, 14, 15,
 16] and it was left open for the critical exponent case $p+1=2(2^*)$ [13; Remark
 3.13].
   The author in [16], proved the
existence of solutions   for $p+1=2(2^*)$ with $N\geq 3$ whenever
the potential function $V(x)$ satisfies some geometry conditions. It
the present paper, we will extend this result for the case $N=2$. It
is well-known that for the semilinear case $(k=0)$,
\begin{eqnarray*}
- \Delta u+V(x)u=g(u), \quad u>0, x \in
  {\BR}^N, \quad \quad \quad
\quad \quad \quad \quad \quad \quad \quad (P)
\end{eqnarray*}

$p+1=2^*$ is the critical exponent when $N \geq 3$. Here is the
definition of the critical growth for  $N=2,$

\begin{itemize}
  \item {\it   Critical growth:  }
 There exists $\beta_0>0$ such that
 $$\lim_{t\rightarrow +\infty} \frac {|g(x,t)|}{\exp(\beta t^2)}=0 \quad \quad \forall \beta >\beta_0,
 \quad  \lim_{t\rightarrow +\infty} \frac {|g(x,t)|}{\exp(\beta t^2)}=+\infty \quad \quad \forall \beta <\beta_0 $$
 \end{itemize}
 uniformly with respect to $x \in \mathbb{R}^2.$
Note that the corresponding critical growth for $N=2$ comes from a
version of Trudinger-Moser inequality in the  whole space ${\BR}^2$
(see [6]) as follows,
\begin{eqnarray*}
\int_{{\BR}^2}\big( \exp (\beta |u|^2)-1 \big )\, dx <+\infty, \quad
\quad \quad \forall u \in H^1({\BR}^2), \beta >0.
\end{eqnarray*}
Also, if $\beta< 4\pi$ and $|u|_{L^2({\BR}^2)}\leq C,$ there exists
a constant $C_2=C_2(C,\beta)$ such that
\begin{eqnarray*}
\sup_{|\nabla u|_{L^2({\BR}^2)}\leq 1}\int_{{\BR}^2}\big( \exp
(\beta |u|^2)-1 \big )\, dx <C_2.
\end{eqnarray*}

There are many results about the existence of solutions for the
subcritical, critical  and the supercritical
exponent case for problem ($P$)(e.g. [1, 3, 4, 5, 19, 22]).\\

 In the case $k >0$, for the subcritical case,  the
existence of a nonnegative solution was proved  for $N=1$ by
Poppenberg, Schmitt and Wang in [18] and for $N\geq 2$ by  Liu and
Wang  in [12].
   In [13] Liu and Wang improved
these results by using a change of variables and treating the new
problem in an Orlicz space. The author in [15], using the idea of
the fibrering method, studied this problem in connection with the
corresponding eigenvalue problem for the laplacian $ -\Del u=V(x) u$
and  proved the existence of multiple solutions for  problem (2).
 It is  established in  [11],
the existence of both one-sign and nodal ground states of soliton
type solutions by the Nehari method. They also established some
regularity of the positive solutions.

As it was mentioned above, for the case $k>0$ with $g(u)=|u|^{p-1}u$
and $N \geq 3$, $p+1=2(2^*)= \frac{4N}{N-2}$
 behaves like a critical exponent for problem (2). This is  because of the nonlinearity term  $-\epsilon k(\Del(|u|^{2}))u$.
 Therefore for problem (2), because of the presence of this nonlinearity term,
 the above definition of  Critical
 growth for $N=2$ changes as
   follows:

\begin{itemize}
  \item {\it   \textbf{Critical growth: } }
 If $N=2,$ there exists $\beta_0>0$ such that
 $$\lim_{t\rightarrow +\infty} \frac {|g(x,t)|}{\exp(\beta t^4)}=0 \quad \quad \forall \beta >\beta_0,
 \quad  \lim_{t\rightarrow +\infty} \frac {|g(x,t)|}{\exp(\beta t^4)}=+\infty \quad \quad \forall \beta <\beta_0, $$
\end{itemize}
uniformly with respect to $x \in {\BR}^2$

Here, we shall study problem (2) with $N= 2$ and show the existence
of positive solutions when the function $g$ has the critical growth.
Before to state the main result, we fix the hypotheses on the
potential function  $V$ and the function $g.$ Indeed,
\begin{enumerate}
\item[{\bf H1:}] $V:{\BR}^2\rightarrow
{\BR}$ is a continuous $1-$periodic function satisfying $V(x)\geq
V_0 >0$ for all $x \in {\BR}^2.$
\item[{\bf H2:}] $g:{\BR}^2 \times {\BR}\rightarrow
{\BR}$ is a nonnegative continuous $1-$periodic function satisfying
$g(x,s)=o_1(s)$ near origin uniformly with respect to $x \in
{\BR}^2.$
\item[{\bf H3:}] $g$ has critical growth at $+\infty,$ namely,
$$g(x,s)\leq C(e^{4\pi s^4}-1) \text{  for all  }(x,s) \in {\BR}^2 \times [0, +\infty).$$
\item[{\bf H4:}] The Ambrosetti-Rabinowitz growth condition: There
exists $\theta >4$ such that
\begin{eqnarray*}
0\leq \theta G (x,t)= \theta \int _0^t g(x,s) \, ds \leq t g(x,t),
\quad \quad t>0.
\end{eqnarray*}
\item[{\bf H5:}] For each fixed $x \in {\BR}^2,$ the function $\frac {g(x,s)}{s}$ is increasing with respect to
$s$, for $s>0.$
\item[{\bf H6:}] There are constants $p>2$ and $C_p$ such that
$$g(x,s)\geq C_ps^{p-1} \text{  for all  }(x,s) \in {\BR}^2 \times [0, +\infty),$$

where
\begin{eqnarray}
C_p & > & \Big [ \frac {\theta (p-2)}{p (\theta-4)} \Big ]^{(p-2)/2}
S_p^p,\\
 S_p& :=&\inf_{u \in H_r^1 ({\BR}^2)\setminus \{0\}}\frac { \big
(\int_{{\BR}^2}(| \nabla u|^2+V_1u^2) \, dx+(\int_{{\BR}^2}u^2|
\nabla u|^2\, dx)^{1/2} \big )^{1/2}}{ \big (\int_{{\BR}^2}|u|^p \,
dx \big )^{1/p} },
\end{eqnarray}
and $V_1:= \max_{x \in {\BR}^2} V(x).$
 \end{enumerate}
It follows from Theorem 1.1 in [12] by some obvious changes  that
the infimum in
(4) attains.\\
 Here is our main Theorem.
\begin{Th} Assume Conditions $ H1- H6.$ Then,
 $(2)$ possesses a nontrivial weak  solution $\tilde{u} \in
H^1({\BR}^N)$ with $\tilde{u}^2\in H^1({\BR}^N).$
\end{Th}
This paper is organized as follows. In Section 2, we reformulate
this problem in an appropriate Orlicz space.  Theorem 1.1 is proved
in Section 3.
\section{Reformulation of the problem and preliminaries}
In this section we assume $N\geq2.$ Denote by $H_{r}^1 ({\BR}^{N})$
the space of radially symmetric functions in
\begin{eqnarray*}
H^{1,2} ({\BR}^N)= \left \{ u \in L^{2}({\BR}^N) : \nabla u \in
L^{2}({\BR}^N) \right \}.
\end{eqnarray*}
 Without loss   of
generality, one can assume $k=1$ in problem (2).
 We formally
formulate problem (2) in a variational structure as follows
$$J(u)=\f{1}{2} \int_{{\BR}^N}(1+u^2)| \nabla u|^2dx+\f{1}{2}
\int_{{\BR}^N}V(x)u^2dx- \int_{{\BR}^N} G(x,u)dx.$$ on the space
$$X=\big \{u\in H^{1,2}({{\BR}}^N): \int_{{{\BR}}^N} V(x) u^2dx<\infty \big
\},$$
which is equipped with the following norm,
$$\|u\|_X=\left \{ \int_{{\BR}^N}| \nabla u|^2dx+ \int_{{{\BR}}^N} V(x) u^2dx\right \}^{\frac{1}{2}}.$$
 Liu and Wang in [13] for the subcritical case,  by making a change of variables treated
this problem in an  Orlicz space. Following their work, we consider
this problem for the supercritical exponent case  in the same Orlicz
space. To  convince the reader we briefly recall some of their
notations and results that are useful in the sequel.

First, we make a change of variables as follows,
$$dv=\sqrt{1+u^2}du, \quad v=h(u)=\f{1}{2}u\sqrt{1+u^2}+\f{1}{2} \ln
(u+\sqrt{1+u^2}).$$ Since $h$ is strictly monotone it has a
well-defined inverse function: $u=f(v)$. Note that
$$h(u)\sim \begin{cases}
u, & |u|\ll 1\\
\f{1}{2}u|u|, & |u|\gg 1, \quad h'(u)=\sqrt{1+u^2},
\end{cases}$$
and
$$f(v)\sim \begin{cases}
v & |v|\ll 1 \\
\sqrt{\f{2}{|v|}}v, & |v|\gg 1, \quad
f'(v)=\f{1}{h'(u)}=\f{1}{\sqrt{1+u^2}}= \f{1}{\sqrt{1+f^2(v)}}.
\end{cases}$$
Also, for some $C_0>0$ it holds
$$L(v):=f(v)^2\sim \begin{cases} v^2& |v|\ll 1,\\
2|v|& |v|\gg 1, \quad L(2v)\leq C_0 L(v),
\end{cases}$$
$L(v)$ is convex, $L'(v)=\f{2f(v)}{\sqrt{1+f(v)^2}}$,
$L''(v)=\f{2}{(1+f(v)^2)^2}>0$.

Using this change of variable, we can rewrite the functional $J(u)$
as
$$\bar {J}(v)=\f{1}{2} \int_{{\BR}^N} |\nabla v|^2dx+ \f{1}{2} \int_{{\BR}^N}V(x)
f(v)^2 dx-\int_{{\BR}^N} G(x,f(v))dx.$$ $\bar {J}$ is defined on the
space
$$H^1_L ({\BR}^N)=\big \{v|   v(x)=v(|x|),   \int_{{\BR}^N}|\nabla v|^2dx <\infty, \int_{{\BR}^N}V(x)
L(v)dx<\infty \big \}.$$ We introduced the Orlicz space (e.g.[20])
$$E_L ({\BR}^N)=\big \{v| \int_{{\BR}^N}V(x) L(v)dx<\infty \big \},$$
equipped with the  norm
$$|v|_{E_L ({\BR}^N)}=\inf_{\zeta>0} \zeta \big (1+\int_{{\BR}^N}(V(x) L(\zeta^{-1}v(x))dx\big ),$$
and define the norm of $H^1_L ({\BR}^N)$ by
$$\|v\|_{H^1_L ({\BR}^N)}= |\nabla v|_{L^2({\BR}^N)}+|v|_{E_L ({\BR}^N)}.$$
Here are  some related facts.
\begin{Prop}
 \begin{enumerate}
 \item[{ (i)}]
 $E_L({\BR}^N)$ is a Banach space.
 \item[{(ii)}]  If $v_n\lo v$ in $E_L({\BR}^N)$, then $\int_{{\BR}^N}V(x)|
L(v_n)-L(v)| dx\lo 0$ and  $\int_{{\BR}^N}V(x)| f(v_n)-f(v)|^2dx\lo
0$.
 \item[{(iii)}]  If $v_n\lo v$ a.e. and $\int_{{\BR}^N}V(x) L(v_n)dx\lo
\int_{{\BR}^N}V(x)L(v)dx$, then $v_n\lo v$ in $E_L({\BR}^N)$.
 \item[{(iv)}]  The dual space $E^*_L({\BR}^N)=L^\infty\cap L_V^2=\{w| w\in L^\infty,
\int_{{\BR}^N}V(x)w^2dx<\infty\}$.
 \item[{ (v)}]  If $v\in E_L({\BR}^N)$, then $w=L'(v)=2f(v)f'(v)\in E^*_L({\BR}^N)$, and
$|w|_{E^*_L}=\sup_{|\phi|_{E_L ({\BR}^N)}\leq 1}(w,\phi)\leq C_1
\big (1+\int_{{\BR}^N} V(x)L(v)dx \big )$, where $C_1$ is a constant
independent of $v$.
 \item[{(vi)}] For $N>2$ the map:$v\lo f(v)$ from $H^1_L({\BR}^N)$ into
$L^q({\BR}^N)$ is continuous for $2\leq q\leq 22^*$ and is compact
for $2< q <22^*.$ Also, for $N=2$, this map is continuous for $q>1.$
 \item[{(vii)}] Suppose $0<V_0\leq V(x)<V_1.$ There exists a
 positive constant $C$ such that
 $$\|u^2\|_{H^{1,2}({\BR}^N)}\leq C (\|v\|_{H^1_L({\BR}^N)}+\|v\|^2_{H^1_L({\BR}^N)}),$$
 where $u=f(v).$
 \end{enumerate}
\end{Prop}
\paragraph{\bf Proof.} See Propositions~(2.1) and (2.2) in [13] for the proof of parts $(i)$ to $(vi)$. We
prove part $ \it {(vii)}$. A direct computation shows that
$$\int |\nabla u^2|^2 \, dx=4 \int \frac{f(v)^2 |\nabla v|^2}{1+f(v)^2} \, dx\leq 4\|v\|_{H^1_L({\BR}^N)}^2.$$
Also, from part (vi)  we have
$$\int u^4 \, dx = \int f(v)^4 \, dx \leq C \|v\|_{H^1_L({\BR}^N)}^4.$$
Now the result is deduced from  the above inequalities. $\Box$

 Hence forth, $\int, H^1, H^1_r,  H^1_L, E_L, L^t, |\cdot|_L$  and $\|\cdot\|$
stand for $\int_{{\BR}^2}$, $H^{1,2}({\BR}^2)$,  $ H^1_r({\BR}^2)$,
$ H^1_L({\BR}^2)$, $ E_L({\BR}^2)$, $L^t({\BR}^2),$ $|\cdot|_{E_L
({\BR}^2)}$ and $\|\cdot\|_{H^1_L({\BR}^2)}$ respectively. In the
following we use $C$ to denote any constant that is independent of
the sequences considered.

\section{Proof of Theorem 1.1}

In this section, we combine the arguments used in [1] and [16]  to
prove Theorem 1.1. The following proposition states some properties
of the functional $\bar {J}.$
 \begin{Prop}
 \begin{enumerate}
 \item[{(i)}]  $\bar{J}$ is well-defined on $H^1_L$.
 \item[{(ii)}]  $\bar{J}$ is continuous in $H^1_L$.
 \item[{(iii)}] $\bar{J}$ is Gauteaux-differentiable in $H^1_L$.
 \end{enumerate}
\end{Prop}

\paragraph{\bf Proof.} The proof is similar to the proof of
Proposition~(2.3) in [13] by some obvious changes.$\square$\\

Here, we prove the existence of a critical point for the functional
$\bar J.$
\begin{Th}
 $ \bar J$ has a critical point  in $H^1_L$, that is,
there exists $0 \neq v\in H^1_L$ such that
$$\int \nabla
v.\nabla \phi dx+\int V(x)f(v)f'(v)\phi dx-\int g(x,f(v))f'(v)\phi
dx=0,$$ for every $\phi\in H^1_L$.
\end{Th}

We use the Mountain Pass Theorem (see [2], [19]) to prove Theorem
3.2.  First, let us define the Mountain Pass value,
$$C_0:=\inf_{\ga\in\Ga} \sup_{t\in[0,1]} \bar{J}(\ga(t)),$$
where
$$\Ga=\{\ga\in C([0,1], H^1_L) | \ga(0)=0, \bar {J}(\ga(1))\leq 0, \ga(1)\neq 0\}.$$
The following lemmas are crucial for the proof of Theorem 3.2.
\begin{Lem}
 The functional $\bar{J}$ satisfies the Mountain Pass
Geometry. \end{Lem}
\paragraph{\bf Proof.}
 We need to show that there exists $0\neq v \in~H^1_L$ such that
$\bar{J}(v)\leq 0$.   Let $0\neq u\in C_{0}^\infty({\BR}^2).$  It is
easy to see that ${J}(tu)\leq 0 $ for the large values of $t.$
Consequently
  $\bar {J} (v) <0$ where $v=h(tu).$ $\square$
\begin{Lem} $C_0$\; is positive.
\end{Lem}
\paragraph{\bf Proof.} Set
\begin{equation*}
S_\rho:=\big \{v\in H^1_L| \int |\nabla v|^2dx+\int V(x)
f(v)^2dx=\rho^2 \big \}.
\end{equation*}

It follows from $H2$ that there exists $\del >0$ such that
\begin{eqnarray}
G(x,s)\leq \f {|s|^2}{4 V_1}, \quad \quad |s|\leq \del.
\end{eqnarray}
Also, for $|s|>\del$,it follows from $H3$ and $H4$ that

\begin{eqnarray}
\del^p G(x,s) \leq |s|^p G(x,s)\leq \f {|s|^{p+1}g(x,s)}{\tht}  \leq
C|s|^{p+1}(e^{4\pi s^2}-1).
\end{eqnarray}
By (5) and (6), we get
\begin{eqnarray}
 G(x,s) \leq \f {|s|^2}{4 V_1} +
\f {C}{\del^p}|s|^{p+1}(e^{4\pi s^2}-1).
\end{eqnarray}

Set $u=f(v),$ with  $v \in S_\rho$ and $\rho <1$. By part (vii) of
Proposition 2.1 $$\|u^2\|_{H^{1,2}({\BR}^N)}\leq C
(\|v\|_{H^1_L({\BR}^N)}+\|v\|^2_{H^1_L({\BR}^N)}),$$  hence, if
$\rho$ is sufficiently small,  it follows from Trudinger-Moser
inequality and the above inequality that
\begin{eqnarray*}
\int (e^{4 q \pi  f(v)^4}-1)\, dx=\int (e^{4 q \pi u^4}-1)\, dx \leq
C
\end{eqnarray*}
for every $q>1$ close to one. Therefore,  it follows from the above
inequality and H\"{o}lder inequality that, $(\frac {1}{q} +\frac
{1}{q'}=1 )$
\begin{eqnarray}
\int |f(v)|^{p+1}(e^{4\pi f(v)^4}-1)\, dx & \leq & \Big (\int
|f(v)|^{q'(p+1)} \, dx \Big )^{\f {1}{q'}} \Big (\int (e^{4q\pi
f(v)^4}-1)
 \, dx \Big )^{\f {1}{q}} \nonumber \\
 & \leq & C \Big (\int
|f(v)|^{q'(p+1)} \, dx \Big )^{\f {1}{q'}}\nonumber \\ & \leq & C
\|v\|^{p+1}
\end{eqnarray}
Taking into account (7) and (8) for each $v \in S_\rho$ with $\rho
<<1$, we have
\begin{eqnarray}
\int G(x,f(v)) \, dx \leq \f{1}{4}\rho ^2+ \f {C}{\del^p}\rho^{p+1}
\end{eqnarray}

Considering  (9) and the fact that $v\in S_\rho$, we obtain
\begin{align*}
\bar{J}(v)&=\f{1}{2}\int |\nabla v|^2dx+\f{1}{2} \int V(x)f(v)^2
dx-\int
G(x,f(v))dx\\
&\geq \f{1}{2}\rho^2-\f{1}{4}\rho^2-\f{C}{\del^P} \rho^{p+1} \geq
\f{1}{8}\rho^2,
\end{align*}
when $0<\rho\leq \rho_0\ll 1$ for some $\rho_0.$  Hence, for $v\in
S_\rho$ with $0<\rho\leq \rho_0$ we have
\begin{equation}
\bar{J}(v)\geq \f{1}{8}\rho^2.
\end{equation}
If $\ga(1)=v$ and $\bar {J}(\ga(1))<0$ then it follows from (10)
that
$$\int |\nabla v|^2dx+ \int V(x)f(v)^2dx >\rho_0^2,$$
thereby giving
$$\sup_{t\in [0,1]}\bar{J}(\ga(t)) \geq \sup_{\ga(t)\in S_{\rho_0}}
\bar{J}(\ga(t))\geq \f{1}{8}\rho_0^2.$$ Therefore $C_0\geq
\f{1}{8}\rho_0^2>0$.$\Box$
\begin{Lem} $C_0$ is bounded from above
   by $\f {\theta-4}{2 \theta}.$
\end{Lem}
\paragraph{\bf Proof.} We fix  a positive radial function $\phi \in
H^1_r$ such that

\begin{eqnarray*}
S_p=\frac { \big (\int_{{\BR}^2}(| \nabla \phi|^2+V_1\phi^2) \,
dx+(\int_{{\BR}^2}\phi^2| \nabla \phi|^2\, dx)^{1/2} \big )^{1/2}}{
\big (\int_{{\BR}^2}|\phi|^p \, dx \big )^{1/p} }, \quad \text { and
} \int_{{\BR}^2}\phi^2| \nabla \phi|^2\, dx \leq 1.
\end{eqnarray*}
Set $\ga_1 (t):= h(t \phi).$ It follows from the definition of the
Mountain Pass value that
$$C_0=\inf_{\ga\in\Ga} \sup_{t\in[0,1]}
\bar{J}(\ga(t))\leq \sup_{t\in[0,1]}
\bar{J}(\ga_1(t))=\sup_{t\in[0,1]} \bar{J}(h(t
\phi))=\sup_{t\in[0,1]} {J}(t\phi).$$

Therefore, we obtain
\begin{align*}
C_0 &\leq
\sup_{t\in[0,1]}J(t\phi) \nonumber \\
&\leq\sup_{t\in [0,1]} \f{ t^2}{2} \int (| \nabla
\phi|^2+V_1\phi^2)\, dx + \f{ t^4}{2} \int
 \phi^2| \nabla \phi|^2 \, dx- \int G(x,t\phi)dx  \nonumber\\
&\leq \sup_{t\in [0,1]} \f{ t^2}{2} \int (| \nabla
\phi|^2+V_1\phi^2)\, dx + \f{ t^2}{2} \Big (\int
\phi^2| \nabla \phi|^2\, dx \Big )^{1/2}- \int G(x,t\phi)dx  \nonumber\\
 &\leq \sup_{t\in [0,1]} \f{ t^2}{2} \int (| \nabla
\phi|^2+V_1\phi^2)\, dx + \f{ t^2}{2} \Big (\int
 \phi^2| \nabla \phi|^2\, dx \Big )^{1/2}- t^p C_p \int \phi^p dx  \nonumber\\
 &=\f{(p-2)S_p^{\f{2p}{p-2}}}{2pC_p^{\f{2}{p-2}}}.
 \end{align*}

 Also, it follows from $H6$ that  $\f{(p-2)S_p^{\f{2p}{p-2}}}{2pC_p^{\f{2}{p-2}}}< \f {\theta-4}{2 \theta}$
 which implies $C_0 <\f {\theta-4}{2 \theta}.$ $\Box$

 The Mountain Pass Theorem guaranties the existence
of a $(PS)_{C_0}$ sequence $\{v_n\},$ that is, $\bar{J}(v_n)\lo C_0$
and $\bar{J}'(v_n)\lo 0$. The following lemma states some properties
of this sequence.
\begin{Lem}
 Suppose $\{v_n\}$ is a $(PS)_{C_0}$ sequence. The
following statements hold.
 \begin{enumerate}
 \item[{(i)}] $\{v_n\}$ is bounded in $H^1_L.$
\item[{(ii)}] If $v_n\geq 0$ converges  weakly to $v $  in $ H^1_L$,  then for every nonnegative test  function $\phi \in  H^1_L $ we have
$$\lim_{n\rightarrow+\infty} \braket{\bar{J}'(v_n)}{\phi}= \braket{\bar{J}'(v)}{\phi}.$$
\end{enumerate}
\end{Lem}
\paragraph{\bf Proof.} Since $\{v_n\}$ is  a $(PS)_{C_0}$ sequence, we have
\begin{eqnarray}
\bar{J}(v_n)&=&\f{1}{2} \int |\nabla v_n|^2 dx+\f{1}{2} \int
V(x)f(v_n)^2 dx- \int G(x,f(v_n))dx\nonumber \\&=&C_0+o(1),
\end{eqnarray}
and
\begin{eqnarray}
\braket{\bar{J}'(v_n)}{\phi}&=&\int \nabla v_n.\nabla \phi dx+\int
V(x)f(v_n)f'(v_n)\phi
dx-\int g(x, f(v_n))f'(v_n) \phi dx \nonumber \\
&=&o(\|\phi\|)
\end{eqnarray}
For part $(i),$ pick $\phi=\f{f(v_n)}{f'(v_n)}=
\sqrt{1+f(v_n)^2}f(v_n)$ as a test function. One can easily deduce
that $|\phi|_L\leq C|v_n|_L$ and
$$|\nabla \phi|=(1+\f{f(v_n)^2}{1+f(v_n)^2}) |\nabla v_n| \leq 2|\nabla v_n|,$$
which implies  $\|\phi\|\leq C\|v_n\|$. Substituting $\phi$ in (8),
gives
\begin{align}
\braket{\bar{J}'(v_n)}{
 \f{f(v_n)}{f'(v_n)}}&= \int (1+\f{f(v_n)^2}{1+f(v_n)^2}) |\nabla
v_n|^2dx+ \int V(x)f(v_n)^2dx \nonumber \\
&\quad - \int g(x,f(v_n))f(v_n)dx \nonumber \\
&=o(\|v_n\|).
\end{align}

 Taking into
account (11), (12) and (13), we have
\begin{align*}
C_0+o(1)+o(\|v_n\|) =&\bar{J}(v_n)-\f{1}{\tht}
\braket{\bar{J}'(v_n)}{
\f{f(v_n)}{f'(v_n)}}\\
=&\f{1}{2}\int |\nabla v_n|^2 dx+\f{1}{2} \int V(x)f(v_n)^2dx
 \\
&-\f{1}{\tht} \int (1+\f{f(v_n)^2}{1+f(v_n)^2}) |\nabla v_n|^2 dx -
\f{1}{\tht}
\int V(x)f(v_n)^2dx \\
=&\int(\f{1}{2}-\f{1}{\tht}(1+\f{f(v_n)^2}{1+f(v_n)^2})) |\nabla
v_n|^2dx+
(\f{1}{2}-\f{1}{\tht}) \int V(x) f(v_n)^2dx \\
\geq & (\f{1}{2}-\f{2}{\tht}) \int |\nabla v_n|^2
dx+(\f{1}{2}-\f{1}{\tht}) \int V(x)f(v_n)^2dx\\
 \geq &
\f{(\tht-4)}{2\tht} \int \big ( |\nabla v_n|^2+V(x)f(v_n)^2\big )\,
dx.
\end{align*}
Since $\tht>4$  it follows from the above that $\int |\nabla
v_n|^2dx+ \int V(x) f(v_n)^2dx$ is bounded and indeed,
\begin{eqnarray}
\limsup_{n\rightarrow\infty} \|v_n\|^2:=K\leq \f {2 \tht
C_0}{\tht-4}.
\end{eqnarray}
It proves part $(i).$

For part $(ii)$,  note first  that it follows from (14) and  Lemma
(3.5) that
\begin{eqnarray*}
\limsup_{n\rightarrow\infty} \|v_n\|^2=K <1.
\end{eqnarray*}

From Trudinger-Moser inequality, there exists $\gamma, q>1$
sufficiently close to one that $T_n(x):=e^{4\pi f(v_n)^2}-1$ is
bounded in  $L^q.$ Since, $v_n\rightarrow v$ a.e. in $\mathbb{R}^2$
so $T_n(x)\rightharpoonup T(x)=e^{4\pi f(v)^2}-1$ weakly in $L^q.$
Now for each $\phi \in H^1_L$ since $H_L^1\subseteq L^t$ for $t>1$
(Proposition 2.1 part (vi)) we have
\begin{eqnarray*}
\int T_n(x)\phi \, dx \rightarrow \int T(x)\phi \,dx.
\end{eqnarray*}
Since, $f$ is increasing and $f(0)=0$, hence $f(v_n)\geq 0$ and $
f(v)\geq0.$  Now it follows from $H_3$  that
\begin{eqnarray*}
g(x, f(v_n))f'(v_n)\phi\leq C T_n(x)\phi.
\end{eqnarray*}
Hence,  the dominated convergence theorem implies
\begin{eqnarray}
\int g(x, f(v_n))f'(v_n)\phi \, dx \rightarrow \int g(x,
f(v))f'(v)\phi \,dx.
\end{eqnarray}
 For the second term on the right hand side of (12), we
have
$$V(x) f(v_n) f'(v_n)\phi \leq  V(x) f(v_n)\phi,$$
and since $v_n\rightharpoonup v$ weakly in $H_1^G$, for the right
hand side of the above inequality we have
$$\lim_{n\rightarrow\infty} \int V(x) f(v_n)\phi \, dx = \int V(x)
f(v)\phi \, dx.$$  Hence by the dominated convergence theorem and
the fact that $v_n \rightarrow v$ a.e. we obtain
\begin{equation}
\lim_{n\rightarrow\infty} \int V(x) f(v_n)f'(v_n)\phi \, dx =  \int V(x) f(v)f'(v)\phi \, dx.
\end{equation}

It follows from (12),  (15) and (16) that
$$\lim_{n\rightarrow+\infty} \braket{\bar{J}'(v_n)}{\phi}=\braket{\bar{J}'(v)}{\phi}.$$
It proves part $(ii)$. $\square$

 Here is a version of Lions' results applicable in our setting.

\begin{Lem} Suppose $v_n\rightarrow 0$ in $H^1_L.$  If there exists $R>0$ such that
\begin{eqnarray}
\liminf_{n\rightarrow\infty} \sup_{y \in \mathbb{R}^2}
\int_{B_R(y)}|f(v_n)|^2 \,dx=0,
\end{eqnarray}
then,
\begin{eqnarray*}
\int f(v_n)g(x,f(v_n)) \,dx \rightarrow 0, \text{ as } n \rightarrow
\infty.
\end{eqnarray*}
\end{Lem}
\paragraph{\bf Proof.}  It follows from (17) and Lemma (4.8) in [7] that
\begin{eqnarray}
 f(v_n)\rightarrow
0 \text { in } L^t \text {  for all  } t\in (2, +\infty).
\end{eqnarray}
    Also, since
 $\limsup_{n\rightarrow\infty} \|v_n\|^2=K <1,$ it follows from
 Trudinger-Moser inequality that
\begin{eqnarray*}
\int (e^{4\gamma\pi f(v_n)^2}-1)\,dx\leq C,
\end{eqnarray*}
for $\gamma >1$ sufficiently close to one. Now, by the same argument
to prove the inequality (7), for any $\epsilon >0$ there exist
constants $C_{\epsilon}$ and $q, \gamma >1$ sufficiently close to
one that
\begin{eqnarray*}
\int f(v_n) g(x, f(v_n))\, dx & \leq & \epsilon \int f(v_n)^2 \,dx +
C_{\epsilon}\int f(v_n)(e^{4\gamma\pi f(v_n)^2}-1)\,dx\\
& \leq & \epsilon C+ C_\epsilon \Big (\int |f(v_n)|^{q'} \, dx \Big
)^{\f {1}{q'}} \Big (\int (e^{4\gamma q\pi f(v_n)^2}-1)
 \, dx \Big )^{\f {1}{q}} \\
 & \leq & \epsilon C+ C_\epsilon \Big (\int
|f(v_n)|^{q'} \, dx \Big )^{\f {1}{q'}}.
\end{eqnarray*}
From the above inequality together with (18), we obtain
\begin{eqnarray*}
\int f(v_n)g(x,f(v_n)) \,dx \rightarrow 0, \text{ as } n \rightarrow
\infty.
\end{eqnarray*}
$\Box$
\begin{Lem}
There exist a sequence $(y_n) $ in $\mathbb{R}^2$ and $R,\epsilon
>0$ such that
\begin{eqnarray*}
\liminf_{n\rightarrow\infty} \int_{B_R(y_n)}|f(v_n)|^2
\,dx>\epsilon,
\end{eqnarray*}
\end{Lem}
\paragraph{\bf Proof.} If $\liminf_{n\rightarrow\infty} \sup_{y \in \mathbb{R}^2}
\int_{B_R(y)}|f(v_n)|^2 \,dx=0,$ it follows from Lemma (3.7) that
$$\int f(v_n)g(x,f(v_n)) \,dx \rightarrow 0, \text{ as } n
\rightarrow \infty.$$  This implies $C_0=0 $ and it is a
contradiction by virtue of Lemma (3.4).$\Box$

\paragraph{\bf Proof of Theorem 3.2.} It follows from Lemma (3.8) the existence of  a sequence $(y_n) $ in
$\mathbb{R}^2$ such that the result of Lemma (3.7) holds for some
$R,\epsilon
>0.$ Without loss of generality, one can assume $(y_n)\subset
\mathbb{Z}^2.$ Now let
 $\tilde{v_n}=v_n(x-y_n).$ Since, $V(.),  g(.,s)$ and $G(.,s)$ are
 $1-$periodic we have
\begin{eqnarray*}
\|v_n\|=\|\tilde{v_n}\|,    J(v_n)=J(\tilde{v_n}) \text { and also
} J'(\tilde{v_n})\rightarrow 0.
\end{eqnarray*}

Since $\|\tilde{v_n}\|$ is bounded, there exists $\tilde{v_0}\in
H^1_L$ such that $\tilde{v_n}\rightarrow \tilde{v_0}$ in $H^1_L.$
Now, it follows from Lemma (3.6) that $J'(\tilde{v_0})=0.$ Also, by
Lemma (3.8) we have
\begin{eqnarray*}
\epsilon \leq  \int_{B_R(0)}|f(\tilde{v_n})|^2 \,dx<2
\int_{B_R(0)}|f(\tilde{v_n})-f(\tilde{v_0})|^2 \,dx+2
\int_{B_R(0)}|f(\tilde{v_0})|^2 \,dx,
\end{eqnarray*}

which implies $  \tilde{v_0}\not \equiv 0.$ $\Box$

\paragraph{\bf Proof of Theorem 1.1.} Proof is  a direct consequence of
Theorem (3.2). Indeed, since $\tilde{v_0}\not \equiv 0$ is a
critical point of $\bar J$, it is easily seen that
$\tilde{u}=f(\tilde{v_0})$ is a nontrivial critical point of $J.$
$\Box$

\end{document}